\documentclass[11pt]{amsart}

\newtheorem{theorem}{Theorem}[section]

\theoremstyle{definition}
\newtheorem{definition}[theorem]{Definition}

\theoremstyle{remark}

\numberwithin{equation}{section}

\begin{document}

\newcommand{\spacing}[1]{\renewcommand{\baselinestretch}{#1}\large\normalsize}
\spacing{1.14}

\title{On the Geometry of  some Para-hypercomplex Lie groups}

\author {H. R. Salimi Moghaddam}

\address{Department of Mathematics, Faculty of  Sciences, University of Isfahan, Isfahan,81746-73441-Iran.} \email{salimi.moghaddam@gmail.com and hr.salimi@sci.ui.ac.ir}



\begin{abstract}
In this paper, firstly we study some left invariant Riemannian
metrics on para-hypercomplex 4-dimensional Lie groups. In each Lie
group, the Levi-Civita connection and sectional curvature have
been given explicitly. We also show these spaces have constant
negative scalar curvatures. Then by using left invariant
Riemannian metrics introduced in the first part, we construct some
left invariant Randers metrics of Berwald type. The explicit
formulas for computing flag curvature have been obtained in all
cases. Some of these Finsler Lie groups are of non-positive flag
curvature.
\end{abstract}

\maketitle

\section{\textbf{Introduction}}\label{Intro}
Hypercomplex and para-hypercomplex structures are important in
differential geometry which have many interesting and effective
applications in theoretical physics. For example the background
objects of HKT-geometry are hypercomplex manifolds. These spaces
appear in $N=4$ supersymmetric model (see \cite{Po,GiPaSt}.). Also
para-hypercomplex structures appear as target manifolds of
hypermultiplets in Euclidean theories with
rigid $N=2$ supersymmetry (see \cite{CoMaMoSa}.).\\
Like M. L. Barberis who has classified invariant hypercomplex
structures on a simply-connected 4-dimensional real Lie group
(\cite{Ba1,Ba2}), N. Bla$\check{z}$i$\acute{c}$ and S.
Vukmirovi$\acute{c}$ have classified 4-dimensional real Lie
algebras which admit a para-hypercomplex structure \cite{BlVu}. In
the first part (section 3) of this paper we consider the connected
Lie groups corresponding to some of these Lie algebras and assume
some left invariant Riemannian metrics on these 4-dimensional Lie
groups. Then for the Riemannian manifolds we compute the
Levi-Civita connection, sectional curvature and scalar curvature.
This shows that all the spaces have constant negative scalar
curvature, also some of them have non-positive sectional curvature.\\
In the second part (section 4) of the paper we try to construct
some invariant Finsler metrics on the Lie groups illustrated in
the first part. Finsler manifolds have many applications in
physics and biology (see \cite{AnInMa,As}.), therefore it can be
important to find some Finsler metrics on these manifolds. On the
other hand the manifolds assumed in this paper are Lie groups so
it is interesting to investigate left invariant Finsler metrics.
The study of invariant Finsler metrics on Lie groups and
homogeneous spaces is one of the attractive fields in Finsler
geometry which has been considered in the recent years (see
\cite{DeHo1,DeHo2,EsSa1,EsSa2,Sa1,Sa2,Sa3}.). Invariant structures
can extricate us from the complicated local coordinates
computations appearing in Finsler geometry. In the pervious works
we have studied invariant metrics on 4-dimensional hypercomplex
Lie groups \cite{Sa4,Sa5}. In this article we continue our studies
on 4-dimensional para-hypercomplex Lie groups and by using
Riemannian metrics introduced in the first part, we construct some
left invariant Randers metric of Berwald type . One of the
fundamental quantities, which is the generalization of sectional
curvature in Riemannian geometry, is flag curvature. Computing
flag curvature of Finsler manifolds in general state is vary
complicated. We obtain the explicit formula for computing the flag
curvature of each Finsler manifold constructed in section 4 and
show that some of the spaces admit left invariant Randers metric
of non-positive flag curvature.


\section{\textbf{Preliminaries}}\label{Prelim}
Let $M$ be a smooth manifold and $\{J_i\}_{i=1,2,3}$ be a family
of fiberwise endomorphisms of $TM$ such that
\begin{eqnarray}
  J_1^2 &=& -Id_{TM}, \\
  J_2^2 &=& Id_{TM}  \ \ \ \ J_2\neq\pm Id_{TM} \\
  J_1J_2&=& -J_2J_1=J_3, \label{JJ}
\end{eqnarray}
and
\begin{eqnarray}
  N_i &=& 0, \ \ \ \ \ \ \ i=1,2,3,
\end{eqnarray}
where $N_i$ is the Nijenhuis tensor corresponding to $J_i$ defined
as follows:
\begin{eqnarray}
   N_1(X,Y)=[J_1X,J_1Y]-J_1([X,J_1Y]+[J_1X,Y])-[X,Y],
\end{eqnarray}
and
\begin{eqnarray}
   N_i(X,Y)=[J_iX,J_iY]-J_i([X,J_iY]+[J_iX,Y])+[X,Y],  \ \ \ \ \
   i=2,3,
\end{eqnarray}
for all vector fields $X,Y$ on $M$. Then the family
$\{J_i\}_{i=1,2,3}$ is called a para-hypercomplex structure
on $M$.\\
In other word a para-hypercomplex structure on a manifold $M$ is a
family $\{J_i\}_{i=1,2,3}$ such that $J_1$ is a complex structure
and $J_i, i=2,3,$ are two non-trivial integrable product
structures on $M$ satisfying in the relation \ref{JJ}.\\
Suppose that $M=G$ is a Lie group. Then we additionally assume
that the para-hypercomplex structure is left invariant, that is:
\begin{definition}
A para-hypercomplex structure $\{J_i\}_{i=1,2,3}$ on a Lie group
$G$ is said to be left invariant if for any $a\in G$
\begin{eqnarray}
  J_i=Tl_a\circ J_i\circ Tl_{a^{-1}}, \ \ \ \ \ i=1,2,3,
\end{eqnarray}
where $Tl_a$ is the differential function of the left translation
$l_a$.
\end{definition}
Also we can consider left invariant metrics on Lie groups.\\
A Riemannian metric $g$ on a Lie group $G$ is called left
invariant if
\begin{eqnarray}
  g(a)(Y,Z)=g(e)(T_al_{a^{-1}}Y,T_al_{a^{-1}}Z), \ \ \ \ \ \forall
  a\in G, \forall Y,Z\in T_aG,
\end{eqnarray}
where $e$ is the unit element of $G$.\\

Let $\frak{g}$ be the Lie algebra of $G$, then the Levi-Civita
connection of the left invariant Riemannian metric $g$ is defined
by the following formula:
\begin{eqnarray}
 2<\nabla_UV,W>=<[U,V],W>-<[V,W],U>+<[W,U],V>,\label{Levi-Civita}
\end{eqnarray}
for any $U,V,W\in\frak{g}$, where $< , >$ is the inner product
induced by $g$ on $\frak{g}$.\\
A Finsler metric on a manifold $M$ is a non-negative function
$F:TM\longrightarrow\Bbb{R}$ with the following properties:
\begin{enumerate}
    \item $F$ is smooth on the slit tangent bundle
    $TM^0:=TM\setminus\{0\}$,
    \item $F(x,\lambda Y)=\lambda F(x,Y)$ for any $x\in M$,
    $Y\in T_xM$ and $\lambda >0$,
    \item the $n\times n$ Hessian matrix $[g_{ij}]=[\frac{1}{2}\frac{\partial^2 F^2}{\partial y^i\partial
    y^j}]$ is positive definite at every point $(x,Y)\in TM^0$.
\end{enumerate}
A special type of Finsler metrics are Randers metrics which have
been introduced by G. Randers in 1941 \cite{Ra}. Randers metrics
are
constructed on Riemannian metrics and vector fields (1-forms).\\
Let $g$ and $X$ be a Riemannian metric and a vector field on a
manifold $M$ respectively such that $\|X\|=\sqrt{g(X,X)}<1$. Then
a Randers metric $F$ can be defined by $g$ and $X$ as follows:
\begin{eqnarray}\label{Randers}
     F(x,Y)=\sqrt{g(x)(Y,Y)}+g(x)(X(x),Y), \ \ \ \ \forall x\in M,
     Y\in T_xM.\label{IRM}
\end{eqnarray}
Similar to the Riemannian case, a Finsler metric $F$ on a Lie
group $G$ is called left invariant if
\begin{eqnarray}
  F(a,Y)=F(e,T_al_{a^{-1}}Y)  \ \ \ \ \ \ \ \forall a \in G, Y\in T_aG.
\end{eqnarray}
We can use left invariant Riemannian metrics and left invariant
vector fields for constructing left invariant Randers metrics on
Lie groups.\\
Suppose that $G$ is Lie group, $g$ is a left invariant Riemannian
metric and $X$ is a left invariant vector field  on $G$ such that
$\sqrt{g(X,X)}<1$. Then we can define a left invariant Riemannian
metric on $G$ by using formula \ref{IRM}. A special family of
Randers metrics (or in general case Finsler metrics) is the family
of Berwaldian Randers metrics. A Randers metric of the form
\ref{IRM} is of Berwald type if and only if the vector field $X$
is parallel with respect to the Levi-Civita connection of $g$. In
these metrics the Chern connection of the Randers metric $F$
coincide with the Levi-Civita connection of the Riemannian metric $g$.\\
One of the important quantities which associates with a Riemannian
manifold is sectional curvature which is defined by the following
formula:
\begin{eqnarray}\label{sectional}
  K(U,Y)=K(P)=\frac{g(R(U,Y)Y,U)}{g(Y,Y).g(U,U)-g^2(Y,U)},
\end{eqnarray}
where $P$ is the 2-plan spanned by $U$ and $Y$ in the tangent
space. Sectional curvature defines another curvature named scalar
curvature as follows:
\begin{eqnarray}\label{scalar}
    S(x)=\sum_{j,k=1, j\neq k}^nK(e_j,e_k),
\end{eqnarray}
where $\{e_1,\cdots,e_n\}$ is an orthonormal basis for $T_xM$ with
respect to the Riemannian metric.\\
The concept of sectional curvature is developed to Finsler
manifolds as follows:
\begin{eqnarray}\label{flag}
  K(P,Y)=\frac{g_Y(R(U,Y)Y,U)}{g_Y(Y,Y).g_Y(U,U)-g_Y^2(Y,U)},
\end{eqnarray}
where $g_Y(U,V)=\frac{1}{2}\frac{\partial^2}{\partial s\partial
t}(F^2(Y+sU+tV))|_{s=t=0}$, $P=span\{U,Y\}$,
$R(U,Y)Y=\nabla_U\nabla_YY-\nabla_Y\nabla_UY-\nabla_{[U,Y]}Y$ and
$\nabla$ is the Chern connection induced by $F$ (see \cite{BaChSh}
and \cite{Sh1}.). This generalization of sectional curvature to
Finsler manifolds is named flag curvature.

\section{\textbf{Left invariant Riamnnian metrics on 4-dimensional para-hypercomlex Lie groups}}

N. Bla\v{z}i\'{c} and S. Vukmirovi\'{c} have studied 4-dimensional
Lie algebras with para-hypercomplex structures. They classified
these spaces in \cite{BlVu}. In this section by using some of
these Lie algebras, we construct some Riemannian Lie groups with
some curvature properties. In the next section by using these
Riemannian metrics, we give some Randers spaces of non-positive
flag curvature.\\
In each case let $G_i$ be the connected 4-dimensional Lie group
corresponding to the considered Lie algebra $\frak{g}_i$ and $<
,>$ is an inner product on $\frak{g}_i$ such that $\{X,Y,Z,W\}$ is
an orthonormal basis for $\frak{g}_i$. Also we use $g$ for the
left invariant Riemannian metric on any $G_i$ induced by $< ,
>$.\\
In any case we suppose that $U=aX+bY+cZ+dW$ and
$V=\tilde{a}X+\tilde{b}Y+\tilde{c}Z+\tilde{d}W$ are any two
independent vectors in $\frak{g}_i$.

\textbf{Case 1.} Let $\frak{g}_1$ be the Lie algebra spanned by
the basis $\{X,Y,Z,W\}$ with the following Lie algebra structure:
\begin{eqnarray}
     [X,Y]=Y  \ \ \ \ , \ \ \ \ [X,W]=W.
\end{eqnarray}
(In any case when we do not write a commutator between the elements of the basis, the commutator is zero.)\\
Now consider the Riemannian manifold $(G_1,g)$. By using formula
\ref{Levi-Civita} for the Levi-Civita connection of $g$ we have:
\begin{eqnarray}
  &&\nabla_XX = 0 \ \ , \ \ \nabla_XY = 0, \ \ \nabla_XZ = 0, \ \ \nabla_XW = 0, \nonumber\\
  &&\nabla_YX = -Y \ \ , \ \ \nabla_YY = X, \ \ \nabla_YZ = 0, \ \ \nabla_YW = 0,\\
  &&\nabla_ZX = 0 \ \ , \ \ \nabla_ZY = 0, \ \ \nabla_ZZ = 0, \ \ \nabla_ZW = 0,\nonumber\\
  &&\nabla_WX = -W \ \ , \ \ \nabla_WY = 0 \ \ , \ \ \nabla_WZ = 0, \ \ \nabla_WW = X\nonumber.
\end{eqnarray}
A direct computation for the curvature tensor shows that:
\begin{eqnarray}
  &&R(X,Y)Y = R(X,W)W = -X,\nonumber\\
  &&R(X,Y)X = -R(Y,W)W = Y,\\
  &&R(X,W)X = R(Y,W)Y = W\nonumber.
\end{eqnarray}
and in other cases $R=0$. Therefore for $U$ and $V$ we have:
\begin{eqnarray}
  R(V,U)U&=&-\{(a\tilde{b}-b\tilde{a})(aY-bX)+(a\tilde{d}-d\tilde{a})(aW-dX)\nonumber\\
  &&+(b\tilde{d}-d\tilde{b})(bW-dY)\}.
\end{eqnarray}
The last equation shows that the sectional curvature is obtained
with the formula:
\begin{eqnarray}
  K(U,V)=-\{(a\tilde{b}-b\tilde{a})^2+(a\tilde{d}-d\tilde{a})^2+(b\tilde{d}-d\tilde{b})^2\}
  \leq 0.
\end{eqnarray}
Also for the scalar curvature $S_1$ we have:
\begin{eqnarray}
  S_1(a)=-6 \ \ \ \ \ \ \ \ \forall a\in G_1.
\end{eqnarray}

\textbf{Case 2.} The Lie algebra of case 2 is of the following
form:
\begin{eqnarray}
     [X,Y]=Z.
\end{eqnarray}
Therefore for the Levi-Civita connection of the Riemannian
manifold $(G_2,g)$ we have:
\begin{eqnarray}
  &&\nabla_XX = 0 \ \ , \ \ \nabla_XY = \frac{1}{2}Z, \ \ \nabla_XZ = -\frac{1}{2}Y, \ \ \nabla_XW = 0, \nonumber\\
  &&\nabla_YX = -\frac{1}{2}Z \ \ , \ \ \nabla_YY = 0, \ \ \nabla_YZ = \frac{1}{2}X, \ \ \nabla_YW = 0,\\
  &&\nabla_ZX = -\frac{1}{2}Y \ \ , \ \ \nabla_ZY = \frac{1}{2}X, \ \ \nabla_ZZ = 0, \ \ \nabla_ZW = 0,\nonumber\\
  &&\nabla_WX = 0 \ \ , \ \ \nabla_WY = 0 \ \ , \ \ \nabla_WZ = 0, \ \ \nabla_WW = 0\nonumber.
\end{eqnarray}
Hence for the curvature tensor of this Levi-Civita connection we
have:
\begin{eqnarray}
  &&-\frac{4}{3}R(X,Y)Y = 4R(X,Z)Z = X,\nonumber\\
  &&\frac{4}{3}R(X,Y)X = 4R(Y,Z)Z = Y,\\
  &&-4R(X,Z)X = -4R(Y,Z)Y = Z\nonumber,
\end{eqnarray}
and so:
\begin{eqnarray}
  R(V,U)U&=&\frac{3}{4}(a\tilde{b}-b\tilde{a})(bX-aY)+\frac{1}{4}(a\tilde{c}-c\tilde{a})(aZ-cX)\nonumber\\
  &&+\frac{1}{4}(b\tilde{c}-c\tilde{b})(bZ-cY).
\end{eqnarray}
This equation shows that the sectional curvature can be obtained
with the formula:
\begin{eqnarray}
  K(U,V)=\frac{1}{4}((a\tilde{c}-c\tilde{a})^2+(b\tilde{c}-c\tilde{b})^2)-\frac{3}{4}(a\tilde{b}-b\tilde{a})^2.
\end{eqnarray}
Also using formula \ref{scalar} for computing scalar curvature
shows that:
\begin{eqnarray}
  S_2(a)=-\frac{1}{2} \ \ \ \ \ \ \ \ \forall a\in G_2.
\end{eqnarray}

\textbf{Case 3.} The Lie algebra structure of $\frak{g}_3$ is of
the following form:
\begin{eqnarray}
     [X,Y]=X.
\end{eqnarray}
Therefore for $(G_3,g)$ we have:
\begin{eqnarray}
  &&\nabla_XX = -Y \ \ , \ \ \nabla_XY = X, \ \ \nabla_XZ = 0, \ \ \nabla_XW = 0, \nonumber\\
  &&\nabla_YX = 0 \ \ , \ \ \nabla_YY = 0, \ \ \nabla_YZ = 0, \ \ \nabla_YW = 0,\\
  &&\nabla_ZX = 0 \ \ , \ \ \nabla_ZY = 0, \ \ \nabla_ZZ = 0, \ \ \nabla_ZW = 0,\nonumber\\
  &&\nabla_WX = 0 \ \ , \ \ \nabla_WY = 0 \ \ , \ \ \nabla_WZ = 0, \ \ \nabla_WW =
  0\nonumber,
\end{eqnarray}
and
\begin{eqnarray}
    R(X,Y)X = Y , R(X,Y)Y = -X.
\end{eqnarray}
Hence for $U$ and $V$ we have:
\begin{eqnarray}
  R(V,U)U=(a\tilde{b}-b\tilde{a})(bX-aY).
\end{eqnarray}
So, like to case (1) we have another Riemannian Lie group of
non-positive sectional curvature with formula:
\begin{eqnarray}
  K(U,V)=-(a\tilde{b}-b\tilde{a})^2\leq 0.
\end{eqnarray}
The scalar curvature $S_3$ of this manifold is of the form:
\begin{eqnarray}
  S_3(a)=-2 \ \ \ \ \ \ \ \ \forall a\in G_3.
\end{eqnarray}

\textbf{Case 4.} In the Lie algebra structure of case (4) there
are two real parameters $\alpha$ and $\beta$. This Lie algebra has
the following structure:
\begin{eqnarray}
     [X,Z]=X \ \ , \ \ [X,W]=Y \ \ , \ \ [Y,Z]=Y \ \ , \ \
     [Y,W]=\alpha X+\beta Y \ \ , \alpha,\beta\in\Bbb{R}.
\end{eqnarray}
Existence of these parameters makes the computing a little
complicated. In this case we have:
\begin{eqnarray}
  &&\nabla_XX = -Z \ \ , \ \ \nabla_XY = \frac{-(1+\alpha)}{2}W, \ \ \nabla_XZ = X, \ \ \nabla_XW = \frac{1+\alpha}{2}Y, \nonumber\\
  &&\nabla_YX = \frac{-(1+\alpha)}{2}W \ \ , \ \ \nabla_YY = -(Z+\beta W), \ \ \nabla_YZ = Y, \ \ \nabla_YW = \frac{(1+\alpha)}{2}X+\beta Y,\\
  &&\nabla_ZX = 0 \ \ , \ \ \nabla_ZY = 0, \ \ \nabla_ZZ = 0, \ \ \nabla_ZW = 0,\nonumber\\
  &&\nabla_WX = \frac{\alpha-1}{2}Y \ \ , \ \ \nabla_WY = \frac{1-\alpha}{2}X \ \ , \ \ \nabla_WZ = 0, \ \ \nabla_WW =
  0\nonumber,
\end{eqnarray}
and
\begin{eqnarray}
  R(X,Y)X = \frac{4-(1+\alpha)^2}{4}Y \ &,& \ R(X,Y)Y = \frac{(1+\alpha)^2-4}{4}X,\nonumber\\
  R(X,Y)Z = 0 \ &,& \ R(X,Y)W = 0,\nonumber\\
  R(X,Z)X = Z \ &,& \ R(X,Z)Y = \frac{1+\alpha}{2}W,\nonumber\\
  R(X,Z)Z = -X \ &,& \ R(X,Z)W = -\frac{1+\alpha}{2}Y,\\
  R(X,W)X = (\frac{1-\alpha^2}{4}+\frac{1+\alpha}{2})W \ &,& \ R(X,W)Y = \frac{1+\alpha}{2}Z+\beta W,\nonumber\\
  R(X,W)Z = -\frac{1+\alpha}{2}Y \ &,& \ R(X,W)W = \frac{(\alpha+1)(\alpha-3)}{4}X-\beta Y,\nonumber\\
  R(Y,Z)X = \frac{1+\alpha}{2}W \ &,& \ R(Y,Z)Y = Z+\beta W,\nonumber\\
  R(Y,Z)Z = -Y \ &,& \ R(Y,Z)W = -\frac{1+\alpha}{2}X-\beta Y,\nonumber\\
  R(Y,W)X =\frac{1+\alpha}{2}Z+\beta W \ &,& \ R(Y,W)Y = \frac{3\alpha^2+4\beta^2+2\alpha-1}{4}W+\beta Z,\nonumber\\
  R(Y,W)Z = -\frac{1+\alpha}{2}X-\beta Y \ &,& \ R(Y,W)W = \frac{-3\alpha^2-4\beta^2-2\alpha+1}{4}Y-\beta X.\nonumber
\end{eqnarray}
Therefore for $U$ and $V$ we have
\begin{eqnarray}
  R(V,U)U=&&-\{(a\tilde{b}-b\tilde{a})(b\frac{(1+\alpha)^2-4}{4}X+a\frac{4-(1+\alpha)^2}{4}Y)\nonumber\\
  &&+(a\tilde{c}-c\tilde{a})(aZ+b\frac{1+\alpha}{2}W-cX-d\frac{1+\alpha}{2}Y)+(a\tilde{d}-d\tilde{a})\nonumber\\
  &&(a\frac{-\alpha^2+2\alpha+3}{4}W+b\frac{1+\alpha}{2}Z+b\beta W-c\frac{1+\alpha}{2}Y\nonumber\\
  &&+d\frac{(1+\alpha)(\alpha-3)}{4}X-d\beta Y)\\
  &&+(b\tilde{c}-c\tilde{b})(a\frac{1+\alpha}{2}W+bZ+b\beta W-c Y-d\frac{1+\alpha}{2}X-d\beta Y)\nonumber\\
  &&+(b\tilde{d}-d\tilde{b})(a\frac{1+\alpha}{2}Z+a\beta W+b\beta Z+b\frac{3\alpha^2+4\beta^2+2\alpha-1}{4}W-c\frac{1+\alpha}{2}X\nonumber\\
  &&-c\beta Y-d\frac{3\alpha^2+4\beta^2+2\alpha-1}{4}Y-d\beta X)\}\nonumber.
\end{eqnarray}
The last equation shows that the sectional curvature can be
obtained with the following formula:
\begin{eqnarray}
  K(U,V)=&&-\{(a\tilde{b}-b\tilde{a})(\tilde{a}b\frac{(1+\alpha)^2-4}{4}+\tilde{b}a\frac{4-(1+\alpha)^2}{4})\nonumber\\
  &&+(a\tilde{c}-c\tilde{a})(a\tilde{c}+b\tilde{d}\frac{1+\alpha}{2}-c\tilde{a}-d\tilde{b}\frac{1+\alpha}{2})\}\nonumber\\
  &&+(a\tilde{d}-d\tilde{a})(a\tilde{d}\frac{-\alpha^2+2\alpha+3}{4}+b\tilde{c}\frac{1+\alpha}{2}+b\tilde{d}\beta-c\tilde{b}\frac{1+\alpha}{2}\nonumber\\
  &&+d\tilde{a}\frac{(1+\alpha)(\alpha-3)}{4}-d\tilde{b}\beta)\\
  &&+(b\tilde{c}-c\tilde{b})(a\tilde{d}\frac{1+\alpha}{2}+b\tilde{c}+b\tilde{d}\beta-c\tilde{b}-d\tilde{a}\frac{1+\alpha}{2}-d\tilde{b}\beta)\nonumber\\
  &&+(b\tilde{d}-d\tilde{b})(a\tilde{c}\frac{1+\alpha}{2}+a\tilde{d}\beta+b\tilde{c}\beta+b\tilde{d}\frac{3\alpha^2+4\beta^2+2\alpha-1}{4}-c\tilde{a}\frac{1+\alpha}{2}\nonumber\\
  &&-c\tilde{b}\beta-d\tilde{b}\frac{3\alpha^2+4\beta^2+2\alpha-1}{4}-d\tilde{a}\beta)\}\nonumber.
  \end{eqnarray}
Also for the scalar curvature $S_4$ we have:
\begin{eqnarray}
  S_4(a)=-\frac{(1+\alpha)^2}{2}-2\beta^2-6 <0\ \ \ \ \ \ \ \ \forall a\in G_4.
\end{eqnarray}

\textbf{Case 5.} The Lie algebra structure of $\frak{g}_5$ is:
\begin{eqnarray}
     [X,Z]=X  \ \ \ \ , \ \ \ \ [Y,W]=Y.
\end{eqnarray}
The Levi-Civita connection of $(G_5,g)$ is as follows:
\begin{eqnarray}
  &&\nabla_XX = -Z \ \ , \ \ \nabla_XY = 0, \ \ \nabla_XZ = X, \ \ \nabla_XW = 0, \nonumber\\
  &&\nabla_YX = 0 \ \ , \ \ \nabla_YY = -W, \ \ \nabla_YZ = 0, \ \ \nabla_YW = Y,\\
  &&\nabla_ZX = 0 \ \ , \ \ \nabla_ZY = 0, \ \ \nabla_ZZ = 0, \ \ \nabla_ZW = 0,\nonumber\\
  &&\nabla_WX = 0 \ \ , \ \ \nabla_WY = 0 \ \ , \ \ \nabla_WZ = 0, \ \ \nabla_WW = 0\nonumber.
\end{eqnarray}
A simple computation for the curvature tensor shows that:
\begin{eqnarray}
  R(X,Z)X = Z \ \ &,& \ \ R(X,Z)Z = -X,\nonumber\\
  R(Y,W)Y = W \ \ &,& \ \ R(Y,W)W = -Y,
  \end{eqnarray}
and
\begin{eqnarray}
  R(V,U)U=-\{(a\tilde{c}-c\tilde{a})(-cX+aZ)+(b\tilde{d}-d\tilde{b})(-dY+bW)\}.
\end{eqnarray}
The above equation shows that the formula for computing sectional
curvature is of the form:
\begin{eqnarray}
  K(U,V)=-\{(a\tilde{c}-c\tilde{a})^2+(b\tilde{d}-d\tilde{b})^2\}
  \leq 0.
\end{eqnarray}
Therefore $(G_5,g)$ is of non-positive sectional curvature and
constant negative scalar curvature $S_5$:
\begin{eqnarray}
  S_5(a)=-4 \ \ \ \ \ \ \ \ \forall a\in G_5.
\end{eqnarray}

\textbf{Case 6.} The last Lie algebra is $\frak{g}_6$ with the
following Lie algebra structure:
\begin{eqnarray}
     [X,Y]=W \ \ , \ \ [X,W]=-Y \ \ , \ \ [Y,W]=-X.
\end{eqnarray}
By using formula \ref{Levi-Civita} for the Levi-Civita connection
of $g$ we have:
\begin{eqnarray}
  &&\nabla_XX = 0 \ \ , \ \ \nabla_XY = \frac{3}{2}W, \ \ \nabla_XZ = 0, \ \ \nabla_XW = -\frac{3}{2}Y, \nonumber\\
  &&\nabla_YX = \frac{1}{2}W \ \ , \ \ \nabla_YY = 0, \ \ \nabla_YZ = 0, \ \ \nabla_YW = -\frac{1}{2}X,\\
  &&\nabla_ZX = 0 \ \ , \ \ \nabla_ZY = 0, \ \ \nabla_ZZ = 0, \ \ \nabla_ZW = 0,\nonumber\\
  &&\nabla_WX = -\frac{1}{2}Y \ \ , \ \ \nabla_WY = \frac{1}{2}X \ \ , \ \ \nabla_WZ = 0, \ \ \nabla_WW = 0\nonumber.
\end{eqnarray}
A computation for the curvature tensor shows that:
\begin{eqnarray}
  &&4R(X,Y)Y = 4R(X,W)W = X,\nonumber\\
  &&-4R(X,Y)X = -\frac{4}{7}R(Y,W)W = Y,\\
  &&-4R(X,W)X = \frac{4}{7}R(Y,W)Y = W\nonumber.
\end{eqnarray}
So we have
\begin{eqnarray}
  R(V,U)U&=&-\frac{1}{4}\{(a\tilde{b}-b\tilde{a})(bX-aY)+(a\tilde{d}-d\tilde{a})(dX-aW)\nonumber\\
  &&+7(b\tilde{d}-d\tilde{b})(-dY+bW)\},
\end{eqnarray}
and
\begin{eqnarray}
  K(U,V)=\frac{1}{4}\{(a\tilde{b}-b\tilde{a})^2+(a\tilde{d}-d\tilde{a})^2-7(b\tilde{d}-d\tilde{b})^2\}.
\end{eqnarray}
The scalar curvature of $(G_6,g)$ is
\begin{eqnarray}
  S_6(a)=-\frac{7}{2} \ \ \ \ \ \ \ \ \forall a\in G_6.
\end{eqnarray}

\section{\textbf{Left invariant Randers metrics on 4-dimensional para-hypercomplex Lie groups}}
In this section we try to use the cases (1) to (6) to construct
invariant Randers metrics of Berwald type. Also in possible cases
we give the explicit formula for computing flag curvature of these
metrics.

\textbf{Case 1.} A direct computation shows that a left invariant
vector field $Q$ is parallel with respect to the Levi-Civita
connection if and only if $Q=qZ$. Now let $0<\|Q\|<1$ because we
would like to construct invariant Randers metric. This shows that
$0<|q|<1$. So the Randers metric defined by $g$ and $Q$ on $G_1$
is of Berwald type. For $g_{\tilde{Y}}(U,V)$ we have (see
\cite{EsSa1}.)
\begin{eqnarray}\label{gtildeY}
  g_{\tilde{Y}}(U,V)&=&g(U,V)+g(X,U).g(X,V)-
  \frac{g(X,\tilde{Y}).g(\tilde{Y},V).g(\tilde{Y},U)}{g(\tilde{Y},\tilde{Y})^{\frac{3}{2}}}\nonumber\\
  &+&\frac{1}{\sqrt{g(\tilde{Y},\tilde{Y})}}\{g(X,U).g(\tilde{Y},V)+g(X,\tilde{Y}).g(U,V)+g(X,V).g(\tilde{Y},U)\}.
\end{eqnarray}
Equation \ref{gtildeY} shows that
\begin{eqnarray}
  g_U(R(V,U)U,V)&=& -(1+qc)\{(a\tilde{b}-b\tilde{a})^2+(a\tilde{d}-d\tilde{a})^2+(b\tilde{d}-d\tilde{b})^2\} \\
  g_U(U,U)&=&(1+qc)^2 \\
  g_U(V,V)&=&1+(q\tilde{c})^2+qc \\
  g_U(U,V)&=&q\tilde{c}(1+qc).
\end{eqnarray}
Now let $P=span\{U,V\}$, therefore for the flag curvature we have
\begin{eqnarray}
  K(P,U) &=&
  \frac{-\{(a\tilde{b}-b\tilde{a})^2+(a\tilde{d}-d\tilde{a})^2+(b\tilde{d}-d\tilde{b})^2\}}{(1+qc)^2}\leq
  0.
\end{eqnarray}
Therefore in the case 1, $(G_1,F)$ is of non-positive flag
curvature.\\

\textbf{Case 2.} Let $\nabla Q=0$ for a left invariant vector
field $Q$. A simple computation shows that $Q=qW$. Also the
assumption $0<\|Q\|<1$ forces us to let $0<|q|<1$. Therefore we
have a left invariant Randers metric of Berwald type which is
constructed by $g$ and $Q$ like formula \ref{Randers}. Also w
have:
\begin{eqnarray}
  g_U(R(V,U)U,V)&=& \{\frac{1}{4}((a\tilde{c}-c\tilde{a})^2+(b\tilde{c}-c\tilde{b})^2)-\frac{3}{4}(a\tilde{b}-b\tilde{a})^2\}(1+qd) \\
  g_U(U,U)&=&(1+qd)^2 \\
  g_U(V,V)&=&1+(q\tilde{d})^2+qd \\
  g_U(U,V)&=&q\tilde{d}(1+qd).
\end{eqnarray}
Therefore the flag curvature formula of this Randers metric is of
the form:
\begin{eqnarray}
  K(P,U) &=&
  \frac{(a\tilde{c}-c\tilde{a})^2+(b\tilde{c}-c\tilde{b})^2-3(a\tilde{b}-b\tilde{a})^2}{4(1+qd)^2}.
\end{eqnarray}

\textbf{Case 3.} A computation for finding left invariant parallel
vector fields with respect to the Levi-Civita connection shows
that the only vector fields which are of the form $Q=q_1Z+q_2W$
are parallel. We consider the vector fields $Q=q_1Z+q_2W$ such
that $0<q_1^2+q_2^2<1$, because we need the condition $0<\|Q\|<1$.
In this case we have:
\begin{eqnarray}
  g_U(R(V,U)U,V)&=& -((a\tilde{b}-b\tilde{a})^2(1+cq_1+dq_2) \\
  g_U(U,U)&=&(1+cq_1+dq_2)^2 \\
  g_U(V,V)&=&1+(\tilde{c}q_1+\tilde{d}q_2)^2+(cq_1+dq_2) \\
  g_U(U,V)&=&(\tilde{c}q_1+\tilde{d}q_2)(1+cq_1+dq_2).
\end{eqnarray}
Now for $P=span\{U,V\}$, the flag curvature is as follows:
\begin{eqnarray}
  K(P,U) &=&
  \frac{-(a\tilde{b}-b\tilde{a})^2}{(1+cq_1+dq_2)^2}\leq 0.
\end{eqnarray}
So in this case $(G_3,F)$ is of non-positive flag curvature.

\textbf{Case 4.} A simple computation shows that $(G_4,g)$ admits
a parallel left invariant vector field $Q$ if and only if
$\alpha=-1$ and $\beta=0$. In this situation $Q=qW$ is the only
family of left invariant vector fields which are parallel with
respect to the Levi-Civita connection of $(G_4,g)$. Let $0<|q|<1$
because of the condition $0<\|Q\|<1$. So we have:
\begin{eqnarray}
  g_U(R(V,U)U,V) &=& -(1+dq)\{(a\tilde{b}-b\tilde{a})^2+(a\tilde{c}-c\tilde{a})^2+(b\tilde{c}-c\tilde{b})^2\}\\
  g_U(U,U) &=& (1+dq)^2 \\
  g_U(V,V)&=& 1+dq+(\tilde{d}q)^2 \\
  g_U(U,V) &=& \tilde{d}q(1+dq).
\end{eqnarray}
Hence for the flag curvature we have:
\begin{eqnarray}
  K(P,U) &=&
  \frac{-\{(a\tilde{b}-b\tilde{a})^2+(a\tilde{c}-c\tilde{a})^2+(b\tilde{c}-c\tilde{b})^2\}}{(1+dq)^2}\leq
  0,
\end{eqnarray}
This case is exactly the case 2 of \cite{Sa5}.

\textbf{Case 5.} In this case a simple computation shows that the
Levi-Civita connection does not admit a left invariant vector
parallel field. Therefore we can not construct Randers metrics as
above cases.

\textbf{Case 6.} The Riemannian Lie group $(G_6,g)$ admits
parallel left invariant vector fields of the form $Q=qZ$. Let
$0<|q|<1$ because we need $0<\|Q\|<1$. Therefore we have:
\begin{eqnarray}
  g_U(R(V,U)U,V)&=& \frac{1+qc}{4}\{(a\tilde{b}-b\tilde{a})^2+(a\tilde{d}-d\tilde{a})^2-7(b\tilde{d}-d\tilde{b})^2\}\\
  g_U(U,U)&=&(1+qc)^2 \\
  g_U(V,V)&=&1+(q\tilde{c})^2+qc \\
  g_U(U,V)&=&q\tilde{c}(1+qc),
\end{eqnarray}
and so
\begin{eqnarray}
  K(P,U) &=&
  \frac{(a\tilde{b}-b\tilde{a})^2+(a\tilde{d}-d\tilde{a})^2-7(b\tilde{d}-d\tilde{b})^2}{4(1+qc)^2}.
\end{eqnarray}
The last equation shows that the flag curvature of $(G_6,F)$ has
not constant sign.

\bibliographystyle{amsplain}

\end{document}